\documentclass[12pt]{article}

\usepackage{amsmath,amsthm,amsfonts}

\parskip=\medskipamount
\def\cinfty#1{C^{\scriptscriptstyle\infty}(#1)}
\def\vectorfields#1{{\cal X}(#1)}
\def\ov#1{\overline{#1}}

\def\fpd#1#2{\frac{\partial #1}{\partial #2}}
\def\R{{\rm I\kern-.20em R}}
\def\sode{{\sc Sode}}

\newtheorem{dfn}{\bf Definition}

\def\baseX#1{{\mathcal X}_{#1}}

\newcommand{\mybox}[1]{\makebox(0,0){\footnotesize{#1}}}
\newlength{\savelen}

\def\Sec{{\mathrm{Sec}}}

\def\tM{\tau_{\scriptscriptstyle M}}
\def\tE{\tau_{\scriptscriptstyle E}}

\def\sV{{\sf V}}
\def\sW{{\sf W}}
\def\sv{{\sf v}}

\def\ss{{\sf s}}

\def\sr{{\sf r}}

\def\se{{\sf e}}

\def\la{{\mathfrak g}}
\def\ua{{\mathfrak u}}
\def\sa{{\mathfrak s}}

\def\pG{\pi_{\scriptscriptstyle G}}

\def\F{{\rm I\kern-.20em F}}

\def\smV{{\scriptscriptstyle V}}

\def\w{\omega}

\begin{document}

\title{A Lie algebroid framework for non-holonomic systems}
\author{Tom Mestdag and Bavo Langerock\thanks{Postdoctoral Fellow of the Fund for Scientific Research - Flanders (Belgium).} \\
{\small Department of Mathematical Physics and Astronomy }\\
{\small Ghent University, Krijgslaan 281, B-9000 Ghent, Belgium}\\
{\small email: \{Tom.Mestdag , Bavo.Langerock\}@UGent.be}}
\date{}

\maketitle

\begin{quote}
{\bf Abstract.} {\small  In order to obtain a framework in which
both non-holonomic mechanical systems and non-holonomic mechanical
systems with symmetry can be described, we introduce in this paper
the notion of a Lagrangian system on a subbundle of a Lie
algebroid.}
\end{quote}



\section{Introduction}
It is well-known that the geometric description of the
Euler-Lagrange equations of a mechanical system, derived from
Hamilton's principle, heavily relies on the intrinsic geometry of
the tangent bundle $TQ$ of the configuration space $Q$. In the
case that the Lagrangian $L$ is invariant under the action of a
Lie group $G$, Hamilton's principle can be reformulated into a
reduced variational principle that takes into account the symmetry
properties of the system. The advantage is that the reduced
variational principle leads to equations defined on a reduced
space, i.e. equations depending on `fewer independent variables'.
Thus, instead of considering the Euler-Lagrange equations on the
total space $Q$, we are interested in the reduced equations, the
so-called Lagrange-Poincar\'e\footnote{We use the terminology of
\cite{CMR,CMR2}.} equations, which are formulated on the Atiyah
quotient bundle $TQ/G\rightarrow Q/G$ (see e.g.\ \cite{Mac}). It
is clear that the geometry of such quotient bundles now becomes of
interest. Weinstein \cite{Wein} has pushed ahead our understanding
in this matter by showing that the geometric structure which lies
at the heart of the Lagrange-Poincar\'e equations is essentially
the same as the one of the Euler-Lagrange equations, namely that
of a Lie algebroid. Therefore, the geometry of `a Lagrangian
system on a Lie algebroid' unifies the geometry of both standard
Lagrangian systems and those where the symmetry properties are
taking into account. In the case of standard Lagrangian systems,
the tangent bundle carries a canonical Lie algebroid structure
which is given by the usual Lie algebra of vector fields on $Q$.
For systems with symmetry, we make use of the so-called Atiyah
algebroid structure on $TQ/G\rightarrow Q/G$  to describe the
evolution equations (see e.g.\ \cite{CMR,DMM}).


In this paper, we mainly deal with mechanical systems that are
subject to some kinematical constraints (i.e. depending on the
velocity), also called non-holonomic constraints. The governing
equations for these models are the so-called Lagrange-d'Alembert
equations. For systems having additional symmetry properties, the
equations of motion can be reduced to the
Lagrange-d'Alembert-Poincar\'e equations (a fairly complete list
of references on non-holonomic mechanics can be found in
\cite{Bloch,Cortes}). The main purpose of this paper is to show
that the Lie algebroid structure constitutes a unifying geometric
structure for describing simultaneously both types of Lagrangian
systems with non-holonomic constraints. Throughout this paper, we
will develop the basic geometric concepts and objects that are
involved.


We believe that at this stage it is instructive to provide a local
version of what will follow (we assume that the reader is familiar
with the concept of a Lie algebroid). We recall the equations of
motion for a non-holonomic mechanical system and take them as the
starting point for further generalisation in the framework of Lie
algebroids, eventually leading to the equations describing a
`non-holonomic Lagrangian system on a Lie algebroid'. Let $L \in
\cinfty{TQ}$ be the Lagrangian of a non-holonomic system. The
constraints are assumed to define a subbundle of the tangent
bundle: $\iota:D\to TQ$, i.e.\ locally we have that
$\iota^k_A(x)w^A = v^k$ where $v^k$ denote the components of
$i(w)$ w.r.t.\ the standard basis $\{\partial/\partial x^k\}$ and
where $w^A$ are the components of $w$ w.r.t.\ some basis of $D$. A
basis of the annihilator space $D^0$ is denoted by $\omega^b
=\omega^b_idx^i$, i.e. $\w^b_i \iota^i_A = 0,\, \forall b,A$. The
Lagrange-d'Alembert principle states that the equations of motion
for such a non-holonomic system are determined by $\delta\int^b_a
L (x^i(t),{\dot x}^i(t))dt =0$, where the variations $\delta x^i$
should satisfy the constraint, i.e. $\delta x \in D_{x(t)}$ for
each $t\in[a,b]$, and moreover $\delta x (a) =\delta x(b) =0$. The
induced equations are called the Lagrange-d'Alembert equations:
\begin{equation} \left\{
\begin{array}{rcl}
\w_i^b \ \dot{x}^i &\!\!\!=\!\!\!& 0, \quad \forall b\\[2mm] \displaystyle
\frac{d}{dt}\left(\fpd{L}{{\dot x}^{i}}\right)- \fpd{L}{x^i}
&\!\!\!=\!\!\!& \displaystyle\lambda_b(t)\omega^b_i(t) \in
D^0_{x(t)},
\end{array}
\right.\label{Leq}
\end{equation}
for some functions $\lambda_b(t)$. It is not hard to see that this
system is equivalent to the system:
\begin{equation} \label{nonholeq}\left\{
\begin{array}{rcl}
\dot{x}^i &\!\!\!=\!\!\!& \iota^i_A(x(t))w^A(t),\\[2mm] \displaystyle
\iota^i_A\frac{d}{dt}\left(\fpd{L}{{\dot x}^{i}}\right)
&\!\!\!=\!\!\!& \iota^i_A \fpd{L}{x^i} , \quad \forall A.
\end{array}
\right. \end{equation} Next, we repeat the above construction for
systems on a Lie algebroid. Assume that a local coordinate chart
$(U,(x,\sv))$ of a Lie algebroid $\tau :\sV\to M$ is given. Given
a Lagrangian $L$ on $\sV$, the Lagrangian equations on the Lie
algebroid are given by:
\begin{equation} \left\{
\begin{array}{rcl}
\dot{x}^i &\!\!\!=\!\!\!& {\rho}^i_{a}(x) {\sv}^{a}, \\[2mm] \displaystyle
\frac{d}{dt}\left(\fpd{L}{{\sv}^{a}}\right) &\!\!\!=\!\!\!&
\displaystyle {\rho}^i_{a} \fpd{L}{x^i} - {C}_{ab}^{c} {\sv}^{b}
\fpd{L}{{\sv}^{c}},
\end{array}
\right.\label{LAeq}
\end{equation}
where $\rho^i_a$ and $C^c_{ab}$ are the structure functions of the
Lie algebroid. We now assume that $\sW$ is a subbundle of $\sV$
with injection $i:\sW\to \sV$. We will use the observation of the
previous paragraph: by contracting the unconstrained Lagrangian
equations with the components of the injection $\iota^a_A$ of the
constraint distribution into the tangent bundle, the non-holonomic
equations can be produced. Assuming that Lagrangian systems
constrained to the subbundle $\sW$ of a Lie algebroid $\sV$ have a
similar behaviour, we can now postulate that the constrained
Lagrangian equations are:
\begin{equation}\label{NHeq} \left\{
\begin{array}{rcl}
\dot x^i(t) &=& {\rho}^i_{a}(x(t)) {\sv}^{a}(t), \\[2mm]
\sv^a(t) &=& i^a_A(x(t)) w^A(t),\\[2mm] \displaystyle
i^a_A\left(\frac{d}{dt}\left(\fpd{L}{{\sv}^{a}}\right) \right)&=&
\displaystyle i^a_A\left({\rho}^i_{a} \fpd{L}{x^i} - {C}_{ab}^{c}
{\sv}^{b} \fpd{L}{{\sv}^{c}}\right).
\end{array} \right.
\end{equation}
We can now write the above equations in terms of derivatives of
the constrained Lagrangian, defined by $L_c(x^i,w^A) =
L(x^i,\sv^a=i^a_A w^A)$. Then, along the solution $\sv=i\circ w$
of (\ref{NHeq}), a straightforward calculation leads to:
\begin{eqnarray*}
\frac{d}{dt}\left(\fpd{L_c}{w^A}\right) &=&
\frac{d}{dt}(i^a_A)\fpd{L}{\sv^a} + i^a_A\frac{d}{dt}
\left(\fpd{L}{\sv^a}\right) \\ &=& \fpd{i^a_A}{x^i} {\dot x}^i
\fpd{L}{\sv^a} + i^a_A \left({\rho}^i_{a} \fpd{L}{x^i} -
{C}_{ab}^{c} {\sv}^{b} \fpd{L}{{\sv}^{c}}\right).
\end{eqnarray*}

If we define $\lambda^i_A=\rho^i_a i^a_A$ then, since
$\displaystyle \partial{L_c}/\partial{x^i} =
\partial{L}/\partial{x^i} + (\partial{L}/\partial{\sv^c}) (\partial{i^c_B}/\partial{x^i})w^B$
, it is clear that the curve $(x(t),w(t))$ is a solution to the
system:
\begin{equation}
\left\{
\begin{array}{rcl}
\dot{x}^i &\!\!\!=\!\!\!& {\lambda}^i_{A}(x) {w}^{A}, \\[2mm] \displaystyle
\frac{d}{dt}\left(\fpd{L_c}{w^A}\right) &\!\!\!=\!\!\!&
\displaystyle \lambda^i_A\fpd{L_c}{x^i} + w^B\left(C^c_{ba}i^b_B
i^a_A  - \lambda^i_A \fpd{i^c_B}{x^i} + \lambda^i_B
\fpd{i^c_A}{x^i}\right) \fpd{L}{\sv^c}.
\end{array}
\right. \label{extension}\end{equation} In the first part of this
paper, we develop all geometric structures required to provide an
intrinsic formulation for this system of equations in
(\ref{extension}). In \cite{Mart}, E.\ Mart\'inez presented a
solid geometrical framework for Weinstein's systems on Lie
algebroids. His approach is very similar to the usual formalism
for Euler-Lagrange equations. The important difference is,
however, that no longer a vector field is the main geometrical
object, but rather a section of a `prolongation bundle' inducing a
unique vector field on the Lie algebroid. In this paper we extend
the framework of Mart\'inez to the above constrained systems
(\ref{extension}). For that purpose, we first define in
Section~\ref{uitwendige} an exterior derivative on a subbundle of
a Lie algebroid. Next, we develop the concept of prolongation
bundles in Section~\ref{sectiedrie}, eventually leading to all
necessary tools for an intrinsic description of (\ref{extension})
in Section~\ref{vergelijkingen}. In our formalism, the system
(\ref{extension}) will be regarded as a section of an appropriate
prolongation bundle. With this section a vector field is
associated, whose integral curves are precisely the solution of
(\ref{extension}). The second part of the paper is devoted to
examples from known dynamical systems that allow a formulation in
terms of the above equations (cf. Section~\ref{examples}). These
examples only deal with autonomous systems with linear
constraints. Finally, we will discuss some of the advantages of
our approach and we will indicate some directions for future work.

\section{Exterior derivatives}\label{uitwendige}

A {\em Lie algebroid} is a vector bundle $\tau: \sV\rightarrow M$
with a real Lie algebra bracket on its set of sections
$[\cdot,\cdot]: \Sec(\tau)\times \Sec(\tau)\rightarrow
\Sec(\tau)$. Moreover there is a linear bundle map
$\rho:{\sV}\rightarrow TM$ (and its natural extension $\rho:
\Sec(\tau)\rightarrow\vectorfields M$) which is related to the
bracket in such a way that, for all ${\ss}, {\sr} \in \Sec(\tau)$,
$f \in \cinfty{M}$
\[ [{\ss} , f {\sr}] = f [{\ss}, {\sr}] + \rho ({\ss})(f) \, {\sr}
\]

is satisfied.

Choose a basis $\{\se_a\}$ for $\Sec(\tau)$ and denote the
corresponding local coordinates on $\sV$ by $(x^i,\sv^a)$. Then,
the structure functions of the Lie algebroid are smooth functions
$\rho^i_a$ and $C^c_{ab}$ on $M$ which satisfy $\rho(\se_a) =
\rho^i_a \fpd{}{x^i}$ and $[\se_a,\se_b] = C^c_{ab} \se_c$. In Lie
algebroid theory, the role of differential forms is played by
sections of exterior powers of the dual vector bundle, i.e.\
skew-symmetric, $\cinfty{M}$-linear maps $\omega: \Sec(\tau)
\times \cdots \times \Sec(\tau)\rightarrow \cinfty{M}$ (with $k$
arguments) will be called {\em $k$-forms on $\Sec(\tau)$} and the
set of all such forms will be denoted by $\bigwedge^k(\tau)$. The
defining properties of a Lie algebroid structure lead to the
definition of an {\em exterior derivative} on $\bigwedge(\tau)$.
Let $\omega\in \bigwedge^k(\tau)$, then the $(k+1)$-form $d\omega$
is given by
\begin{eqnarray}
d\omega({\ss}_1,\ldots,{\ss}_{k+1}) &=&
\sum_{i=1}^{k+1}(-1)^{i-1}\rho({\ss}_i)\Big(
\omega({\ss}_1,\ldots,\hat{{\ss}_i},\ldots,{\ss}_{k+1})\Big) \nonumber \\
&& \mbox{} + \sum_{1\leq i<j\leq k+1}(-1)^{i+j}\omega([
{\ss}_i,{\ss}_j],{\ss}_1,\ldots,\hat{{\ss}_i},
\ldots,\hat{{\ss}_j},\ldots,{\ss}_{k+1}). \label{d}
\end{eqnarray}
The operator $d$ has the property $d(\omega_1 \wedge \omega_2) =
d\omega_1 \wedge \omega_2 + (-1)^{k_1} \omega_1 \wedge d\omega_2$,
with $k_1$ the degree of $\omega_1$, and moreover $d^2=0$.
Locally, $dx^i = \rho^i_a {\se}^a$ and $d{\se}^c =
-\frac{1}{2}C^c_{ab} {\se}^a\wedge{\se}^b$, where the set
$\{\se^a\}$ stands for the basis of $\Sec(\tau^*)$ which is dual
to $\{\se_a\}$. More details on the properties of the exterior
derivative can be found in \cite{DMM,Mac, Mart}.

Suppose now that a vector subbundle $\mu: \sW \rightarrow M$ of
$\tau$ is given. It is obvious that each coordinate system
$(x^i,w^A)$ on $\sW$ can be extended to a coordinate system
$(x^i,\sv^A= w^A, \sv^\alpha)$ on $\sV$. However, in what follows,
it will be more convenient to consider the basis $\{\se_a\}$ of
$\Sec(\tau)$ to be a priori given. Then, the sections of an
arbitrary basis $\{ e_A \}$ of $\Sec(\mu)$, can be written as $e_A
= i^a_A \se_a$ and the injection can be denoted by $i: \sW
\rightarrow \sV; (x^i,w^A) \mapsto (x^i, \sv^a=i^a_A w^A)$. We
will use $\lambda$ for the restriction of $\rho$ to $\sW$:
$\lambda=\rho \circ i : \sW \rightarrow TM; (x^i,w^A)\mapsto
(x^i,{\dot x}^i = \lambda^i_A w^A = \rho^i_a i^a_A w^A)$.

Forms on $\Sec(\tau)$ can be pulled back to forms on $\Sec(\mu)$.
Indeed, if $\omega$ is a $k$-form then $i^*\omega$, defined by
\[
i^*\omega (W_1, \ldots, W_k) = \omega (iW_1, \ldots , iW_k),
\qquad \qquad W_i \in \Sec(\mu)
\]

is a $k$-form on $\Sec(\tau)$. By composing the exterior
derivative $d$ with $i^*$ we can define a mapping from
$\bigwedge^k(\tau)$ to $\bigwedge^{k+1}(\mu)$, which is denoted by
$\delta$. Thus, if $\omega$ is a $k$-form on $\Sec(\tau)$, then
\begin{equation} \delta\omega=i^*d\omega \label{extderdelta},
\end{equation}
is a $(k+1)$-form on $\Sec(\mu)$. The operator $\delta$ satisfies
the rule
\[
\delta(\omega_1 \wedge \omega_2) = \delta\omega_1 \wedge
i^*\omega_2 +(-1)^{k_1} i^*\omega_1 \wedge \delta\omega_2,
\]
and could be called a derivative along $i$ for this reason, but we
will simply refer to $\delta$ as an exterior derivative.
Obviously, $\delta\circ d =0$. In the above introduced
coordinates,
\[
\delta x^i = \rho^i_a i^a_A e^A = \lambda^i_A e^A \qquad
\mbox{and} \qquad \delta \se^a = -\frac{1}{2}D^a_{BC} e^B \wedge
e^C,
\]
where $D^a_{BC} = C^a_{bc} i^b_B i^c_C$.

To illustrate the above notions, we look at the case where $\sW$
is a {\em Lie subalgebroid} of $\sV$. A more general definition of
a Lie subalgebroid $\sW \rightarrow N$ of $\sV \rightarrow M$ can
be found in \cite{Higg}. Here, we will only consider the case that
the base manifolds coincide. A Lie algebroid $\mu:\sW\rightarrow
M$ (with anchor map $\lambda$ and structure functions $D^C_{AB}$)
is a Lie subalgebroid of $\tau:\sV \rightarrow M$ (with anchor map
$\rho$) if there exist a injective morphism $i: \sW \rightarrow
\sV$ such that $\lambda = \rho \circ i$ and
\begin{equation} \label{subalg}
i[W_1,W_2] = [i(W_1),i(W_2)], \qquad\qquad W_i \in \Sec(\mu).
\end{equation}
In coordinates the above expression reads:
\begin{equation} \label{substructure}
i_C^c D^C_{AB} = D^c_{AB} - \lambda^i_A \fpd{i^c_B}{x^i} +
\lambda^i_B \fpd{i^c_A}{x^i}.
\end{equation}
There are now three exterior derivatives around: $d$, $\delta$ and
the exterior derivative $d^\mu$ on the Lie algebroid $\mu$. The
condition (\ref{subalg}) can equivalently be written as
\begin{equation}
d^\mu \circ i^* = i^* \circ d = \delta.
\end{equation}
Indeed, for functions on $M$ it is obviously satisfied, while for
1-forms $\theta\in\bigwedge^1(\tau)$
\[
d^\mu i^*\theta(W_1,W_2) = \lambda(W_1) \theta(i W_2) -
\lambda(W_2) \theta(W_1) - \theta(i[W_1,W_2]).
\]
This is $d\theta(iW_1,iW_2)$ if and only if (\ref{subalg}) is
satisfied. The proof then follows from induction. It is now also
obvious that for Lie subalgebroids
\[ d^\mu \circ \delta= i^* \circ d^2 = (d^\mu)^2 \circ i^*= 0.
\]

\section{Prolongation bundles}\label{sectiedrie}

As stated in the introduction, an important role in our
description will be played by the so-called {\em prolongation
bundle} associated to a Lie algebroid. Sections of such bundles
will eventually lead to the vector field on the Lie algebroid that
generates the equations in~(\ref{extension}). Furthermore, these
bundles allow a generalisation of all intrinsic objects defined on
the tangent bundle needed to write down the standard Lagrange
equations. We start with defining the prolongation of an arbitrary
bundle, and then continue with the introduction of all the
structures that lead to an intrinsic formulation
of~(\ref{extension}).

Let $\pi:E\rightarrow M$ be an arbitrary fibre bundle. The
prolongation of $\pi$ by an anchored vector bundle $\tau:\sV
\rightarrow M$ is the (vector) bundle $\pi^\rho:T^\rho E
\rightarrow E$. Here, the total space $T^\rho{E}$ is the pullback
manifold $\rho^*TE= \{({\sv},X_e) \in {\sV} \times TE \mid \,
\rho(\sv) = T\pi(X_e) \}$. The natural projections of $T^\rho E$
onto $\sV$ and $TE$ will be denoted by $\pi^2$ and $\rho^\pi$,
respectively (see the diagram below). On the other hand, the
bundle projection $\pi^\rho$ is given by $\tE\circ\rho^1$, i.e.\
$\pi^\rho({\sv},X_e)=e$ (see also e.g.\ \cite{DMM,Higg,
Mart,MaMeSa}).

\setlength{\unitlength}{12pt}

\begin{picture}(40,13)(-8,2)

\put(7.5,9.5){\vector(1,0){6.6}} \put(11,13){\vector(1,-1){3}}
\put(7,10.2){\vector(1,1){2.8}}

\put(6.5,9.5){\mybox{$T^\rho{E}$}} \put(14.5,9.5){\mybox{$E$}}
\put(10.5,13.7){\mybox{$TE$}}

\put(7,3.5){\vector(1,0){7}} \put(11,7){\vector(1,-1){3}}
\put(7,4){\vector(1,1){3}}

\put(6.5,3.5){\mybox{${\sV}$}} \put(14.5,3.5){\mybox{$M$}}
\put(10.5,7.5){\mybox{$TM$}}

\put(6.5,8.8){\vector(0,-1){4.7}} \put(10.5,13){\vector(0,-1){5}}
\put(14.5,9){\vector(0,-1){5}}

\put(11,3){\mybox{$\tau$}} \put(8.0,5.7){\mybox{$\rho$}}
\put(13.3,5.7){\mybox{$\tM$}}

\put(11.2,9){\mybox{$\pi^\rho$}} \put(11.4,11){\mybox{$T\pi$}}

\put(8.0,11.8){\mybox{$\rho^\pi$}} \put(13.6,11.5){\mybox{$\tE$}}

\put(6.0,7){\mybox{$\pi^2$}} \put(15.2,7){\mybox{$\pi$}}

\end{picture}
\setlength{\unitlength}{\savelen}

Suppose we fix a bundle adapted coordinate chart
$(U,(x^i,y^\alpha))$ about a point $e\in E$. It is not difficult
to see that $\{ {\mathcal X}_a, {\mathcal Y}_\alpha \}$, with
\begin{equation}\baseX{a} (e)= \left({ \se}_a
(\pi(e)),\left.\rho^i_{a}(x)\fpd{}{x^i} \right|_e\right)\qquad
\mbox{and} \qquad {\mathcal Y}_{\alpha}(e) = \left(
0(\pi(e)),\left.\fpd{}{y^{\alpha}}\right|_e\right), \label{basis}
\end{equation}
is a basis of $\Sec(\pi^\rho)$, induced by the basis $\{{\se}_a\}$
of $\Sec(\tau)$ and the basis $\{\fpd{}{x^i}, \fpd{}{y^\alpha}\}$
of $\vectorfields{E}$. The set ${\mathcal V}^\rho E = \{(0,X_e)\in
T^\rho E\}$ is called the bundle of vertical elements of $T^\rho
E$. A basis for the set of vertical sections $\rm{Ver}(\pi^\rho)$
is given by the $\{ {\mathcal Y}_\alpha \}$. In what follows, we
will use $\{ {\mathcal X}^a, {\mathcal Y}^\alpha \}$ for the basis
of $\Sec((\pi^\rho)^*)$ that is dual to $\{ {\mathcal X}_a,
{\mathcal Y}_\alpha \}$.

Suppose now that $\tau:\sV\rightarrow M$ is a Lie algebroid. The
Lie algebroid structure on $\tau$ can be naturally extended to a
Lie algebroid structure on the prolongation bundle $\pi^\rho:
T^\rho E \rightarrow E$. The anchor map of this Lie algebroid is
$\rho^\pi: T^\rho E \rightarrow TE$, with
\[
\rho^\pi\big(\baseX a\big) = \rho^i_a \fpd{}{x^i} \qquad
\mbox{and} \qquad \rho^\pi\big({\mathcal Y}_\alpha \big)=
\fpd{}{y^\alpha}
\]
and the bracket is given by
\[
[\baseX a,\baseX b] = C^c_{ab}\baseX c, \qquad [\baseX a,
{\mathcal Y}_\beta] = 0, \qquad [{\mathcal Y}_\alpha, {\mathcal
Y}_\beta ]=0
\] (for more details, see e.g.\ \cite{DMM,Mart}).

Let $\mu: \sW \rightarrow M$ be a subbundle of $\tau$. We will
need three different
 prolongation bundles. The first is the one where $\pi$ is
also $\tau:\sV \rightarrow M$, that is
$\tau^\rho:T^\rho\sV\rightarrow \sV$. Since this bundle has a Lie
algebroid structure, it also has an exterior derivative $d$,
locally characterized by the relations
\[
d x^i = \rho^i_a {\mathcal X}^a, \quad d \sv^a = {\mathcal V}^a,
\quad d {\mathcal X}^c = -\frac{1}{2} C^c_{ab} {\mathcal X}^a
{\wedge \mathcal X}^b \quad \mbox{and}\quad d {\mathcal V}^a = 0.
\]
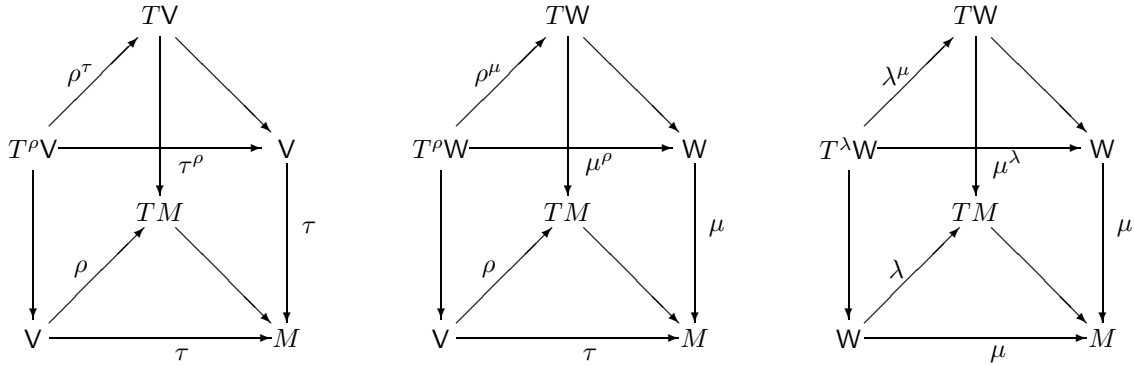
\setlength{\unitlength}{12pt}

\vspace{1cm}

\begin{tabular}{ccc}
\begin{minipage}{5cm}

\begin{center}
\begin{picture}(11,13)

\put(2.3,9.5){\vector(1,0){6.4}} \put(6,13){\vector(1,-1){3}}
\put(2,10.2){\vector(1,1){2.8}}

\put(1.5,9.5){\mybox{$T^\rho{\sV}$}} \put(9.5,9.5){\mybox{$\sV$}}
\put(5.5,13.7){\mybox{$T\sV$}}

\put(2,3.5){\vector(1,0){7}} \put(6,7){\vector(1,-1){3}}
\put(2,4){\vector(1,1){3}}

\put(1.5,3.5){\mybox{$\sV$}} \put(9.5,3.5){\mybox{$M$}}
\put(5.5,7.5){\mybox{$TM$}}

\put(1.5,8.8){\vector(0,-1){4.7}} \put(5.5,13){\vector(0,-1){5}}
\put(9.5,9){\vector(0,-1){5}}

\put(6.2,3){\mybox{$\tau$}} \put(3.0,5.7){\mybox{$\rho$}}

\put(6.5,9){\mybox{$\tau^\rho$}}
\put(3.0,11.8){\mybox{$\rho^\tau$}}

\put(10.2,7){\mybox{$\tau$}}

\end{picture}

\end{center}

\end{minipage}

&

\begin{minipage}{5cm}
\begin{center}

\begin{picture}(11,13)

\put(2.4,9.5){\vector(1,0){6.4}} \put(6,13){\vector(1,-1){3}}
\put(2,10.2){\vector(1,1){2.8}}

\put(1.5,9.5){\mybox{$T^\rho{\sW}$}} \put(9.5,9.5){\mybox{$\sW$}}
\put(5.5,13.7){\mybox{$T\sW$}}

\put(2,3.5){\vector(1,0){7}} \put(6,7){\vector(1,-1){3}}
\put(2,4){\vector(1,1){3}}

\put(1.5,3.5){\mybox{$\sV$}} \put(9.5,3.5){\mybox{$M$}}
\put(5.5,7.5){\mybox{$TM$}}

\put(1.5,8.8){\vector(0,-1){4.7}} \put(5.5,13){\vector(0,-1){5}}
\put(9.5,9){\vector(0,-1){5}}

\put(6.2,3){\mybox{$\tau$}} \put(3.0,5.7){\mybox{$\rho$}}

\put(6.5,9){\mybox{$\mu^\rho$}} \put(3.0,11.8){\mybox{$\rho^\mu$}}

\put(10.2,7){\mybox{$\mu$}}

\end{picture}

\end{center}
\end{minipage}

&

\begin{minipage}{5cm}
\begin{center}

\begin{picture}(11,13)

\put(2.4,9.5){\vector(1,0){6.4}} \put(6,13){\vector(1,-1){3}}
\put(2,10.2){\vector(1,1){2.8}}

\put(1.5,9.5){\mybox{$T^\lambda \sW$}}
\put(9.5,9.5){\mybox{$\sW$}} \put(5.5,13.7){\mybox{$T\sW$}}

\put(2,3.5){\vector(1,0){7}} \put(6,7){\vector(1,-1){3}}
\put(2,4){\vector(1,1){3}}

\put(1.5,3.5){\mybox{$\sW$}} \put(9.5,3.5){\mybox{$M$}}
\put(5.5,7.5){\mybox{$TM$}}

\put(1.5,8.8){\vector(0,-1){4.7}} \put(5.5,13){\vector(0,-1){5}}
\put(9.5,9){\vector(0,-1){5}}

\put(6.2,3){\mybox{$\mu$}} \put(3.0,5.7){\mybox{$\lambda$}}

\put(6.5,9){\mybox{$\mu^\lambda$}}
\put(3.0,11.8){\mybox{$\lambda^\mu$}}

\put(10.2,7){\mybox{$\mu$}}

\end{picture}

\end{center}
\end{minipage}

\end{tabular}

\setlength{\unitlength}{\savelen}

The second prolongation of interest is $\mu^\rho: T^\rho \sW
\rightarrow \sW$. Again, $\mu^\rho$ is a Lie algebroid with a
corresponding exterior derivative, denoted by $\tilde d$. If $\{
{\tilde{\mathcal X}}_a, {\tilde{\mathcal W}}_A \}$ stands for the
basis (\ref{basis}) in this situation, then
\[
\tilde d x^i = \rho^i_a {\tilde{\mathcal X}}^a, \quad \tilde d w^A
= {\tilde{\mathcal W}}^A, \quad \tilde d {\tilde{\mathcal X}}^c =
-\frac{1}{2} C^c_{ab} {\tilde{\mathcal X}}^a \wedge
{\tilde{\mathcal X}}^b \quad \mbox{and}\quad \tilde d
{\tilde{\mathcal W}}^A = 0.
\]
$\mu^\rho$ is in fact a Lie subalgebroid of $\tau^\rho$ (in the
more general sense of \cite{Higg}, see also \cite{DMM}).

Last but not least, we will also need the prolongation
$\mu^\lambda: T^\lambda \sW \rightarrow \sW$, which is a vector
subbundle of the Lie algebroid $T^\rho \sW \rightarrow \sW$. Let
$\{ {\mathcal X}_A, {\mathcal W}_A \}$ stand for the basis of
$\Sec(\mu^\lambda)$, induced by the bases $\{e_A\}$ and
$\{\fpd{}{x^i},\fpd{}{w^A}\}$. Then, the injection is given by the
map $I: T^\lambda \sW \rightarrow T^\rho \sW ; (w_1, X_{w_2})
\mapsto (i(w_1), X_{w_2})$ which, locally, is of the form
\[
I({\mathcal X}_A) = i^a_A {\tilde{\mathcal X}}_a \qquad
\mbox{and}\qquad I({\mathcal W}_A) = {\tilde{\mathcal W}}_A.
\]
In view of the discussion of the previous section, we can use $I$
to introduce an `exterior derivative' $\delta$, defined as in
(\ref{extderdelta}), which will map forms on $\Sec(\mu^\rho)$ to
forms on $\Sec(\mu^\lambda)$. Locally,
\[
\delta x^i = \lambda^i_A {\mathcal X}^A, \quad \delta w^A =
{\mathcal W}^A, \quad \delta {\tilde{\mathcal X}}^c = -\frac{1}{2}
D^c_{AB} {\mathcal X}^A {\wedge \mathcal X}^B \quad
\mbox{and}\quad \delta {\tilde{\mathcal W}}^A = 0.
\]
A second important map $T^\rho i:T^\rho \sW \rightarrow T^\rho
\sV$ is the linear bundle map over $i:\sW\rightarrow\sV$, given by
$T^\rho i(\sv,X_w) = (\sv, T_{i(w)}i(X_w))$.  Let an element $Z\in
T^\rho_w \sW$ be given by $Z = X^a{\tilde{\mathcal X}}_a(w) + W^A
{\tilde{\mathcal W}}_A(w)$, then
\[
 T^\rho i(Z) = X^a {\mathcal X}_a (i(w)) + \big(\rho^i_b X^b \fpd{i^a_A}{x^i}w^A + i^a_A
 W^A\big) {\mathcal V}_a(i(w)).
\]

Both prolongations $\tau^\rho: T^\rho \sV \rightarrow \sV$ and
$\mu^\lambda: T^\lambda \sW \rightarrow \sW$ carry further
interesting canonical objects. For example, there exists a
naturally defined {\em vertical lift} ${}^\smV: \tau^*\sV
\rightarrow {\mathcal V}^\rho \sV \subset T^\rho \sV$. Indeed, for
$({\sf a},\sv)\in \tau^*\sV$, we can define first a vertical
element ${\sf v}_{\sf a}^v \in T_{\sf a}\sV$ by means of its
action on functions $f\in\cinfty{\sV}$,
\begin{equation}\nonumber
{\sv}_{\sf a}^v(f) = \left. \frac{d}{dt} f({\sf a}+t {\sv})
\right|_{t=0}.
\end{equation}
Then, the element ${\sv}_{\sf a}^\smV\in {\mathcal V}^\rho \sV$ is
defined as $({\sf 0},{\sv}_{\sf a}^v)$. We have $\se_a^\smV =
{\mathcal V}_a$. Evidently, there exists a similar notion for the
prolongation $\mu^\lambda$, here the vertical lift is a map
${}^\smV: \mu^*\sW \rightarrow {\mathcal V}^\lambda \sW \subset
T^\lambda \sW$. Next, there are also two so-called {\em vertical
endomorphisms} $S^\tau = {}^\smV \circ j_1: \Sec(\tau^\rho)
\rightarrow Ver(\tau^\rho)$ and $S^\mu = {}^\smV \circ j_2:
\Sec(\mu^\lambda) \rightarrow Ver(\mu^\lambda)$, where $j_1$
stands for the projection $(\sv_2,X_{\sv_1})\in T^\rho\sV\mapsto
(\sv_1,\sv_2)$ and $j_2: T^\lambda \sW \rightarrow \mu^*\sW$ can
be defined analogously. In the above introduced bases, $S^\tau=
{\mathcal X}^a \otimes {\mathcal V}_a$ and $S^\mu = {\mathcal X}^A
\otimes {\mathcal W}_A$. Remark that, although the projection
$j_3:T^\rho\sW \rightarrow \mu^*\sV$ can easily be defined, there
exists no vertical endomorphism on $T^\rho\sW$ because there is no
vertical lift which maps elements in $\mu^*\sW$ to ${\mathcal
V}^\rho \sW$. Two final important objects are the Liouville
sections ${\mathcal C}^\tau= \sv^a {\mathcal V}_a \in
\Sec(\tau^\rho)$ and ${\mathcal C}^\mu = w^A {\mathcal W}_A \in
\Sec(\mu^\lambda)$.

\section{Lagrange equations on a subbundle of a Lie
algebroid}\label{vergelijkingen}

In this section we formulate the equations for a non-holonomic
system on a Lie algebroid. The main purpose is to arrive at an
equation that uniquely determines a vector field on $\sW$ whose
integral curves are precisely the solutions to the constrained
equations \eqref{extension}.  We will only deal with {\em regular}
Lagrangians $L\in\cinfty{\sV}$, i.e.\ we will assume that the
matrix $({\partial^2 L}/{\partial \sv^\alpha\partial \sv^\beta})$
is regular at every point.

Let's recall first briefly Mart\'inez's description of Lagrangian
systems (\ref{LAeq}) on a Lie algebroid. An important subclass of
sections of the prolongation bundle $\tau^\rho: T^\rho\sV
\rightarrow \sV$ are the so-called {\em pseudo-\sode s}
\cite{MaMeSa, MeSa}, or simply \sode s in \cite{DMM,Mart}. They
are sections $\Gamma$ of $\tau^\rho$ such that $\tau^2 \circ
\Gamma = id$ ($\tau^2$ is the projection $(\sv,V)\in T^\rho
\sV\mapsto \sv\in \sV$). Locally, a pseudo-\sode\, is of the form
\[
\Gamma= \sv^a {\mathcal X}_a + f^a {\mathcal V}_a.
\]
It is not difficult to see that pseudo-\sode s are in a one-to-one
correspondence with vector fields $X$ on $\sV$ with the property
that $\rho(\sv)=T\tau(X(\sv))$. Keeping this in mind, by defining
a pseudo-\sode, one can give an intrinsic description of the
equations (\ref{LAeq}). A regular Lagrangian $L$ defines a
function
\[
E_L=\rho^\tau(C^\tau) L - L = \sv^a\fpd{L}{\sv^a}(\sv) -L(\sv)
\]
on $\sV$ and a 1-form
\begin{equation} \theta_L=S^\tau(dL)=
\fpd{L}{\sv^a}(\sv){\mathcal X}^a \label{thetaL}
\end{equation}
on
$\Sec(\tau^\rho)$. The dynamics are then given by a pseudo-\sode\,
$\Gamma$ of the prolongation $\tau^\rho$ that solves the equation
\begin{equation} \label{symplectic}
i_\Gamma d\theta_L = -dE_L.
\end{equation}
Solutions of the equations (\ref{LAeq}) are then nothing but
integral curves of the associated vector field
$\rho^\tau(\Gamma)\in\vectorfields{\sV}$.

For the constrained systems (\ref{extension}), we wish to
preserve, as much as possible, the structure of the equation
(\ref{symplectic}). By analogy, it is clear that we should
represent the dynamics by a section of the prolongation bundle
$\mu^\lambda: T^\lambda \sW \rightarrow \sW$. It is easy to see
that pseudo-\sode s $\Gamma$ on this bundle have locally the form
\[
\Gamma= w^A {\mathcal X}_A + f^A {\mathcal V}_A.
\]
Crucial in (\ref{symplectic}) is the exterior derivative $d$ of
the prolongation $\tau^\rho$. Unfortunately, the prolongation
bundle $\mu^\lambda$ does not carry a Lie algebroid structure and
therefore there is no available exterior derivative. The next best
thing is the above introduced operator $\delta$. As a consequence,
the analogues ${\tilde E}_L$ and ${\tilde\theta}_L$ of the
function $E_L$ and the 1-form $\theta_L$ to the constrained case
should be sought, respectively, among the functions on $\sW$ and
the 1-forms on $\Sec(\mu^\rho)$. For ${\tilde E}_L$ we can simply
take the restriction of $E_L$ to $\sW$. If we define the
constrained Lagrangian $L_c\in \cinfty{\sW}$ as the restriction of
$L$ to $\sW$, $L_c(w) = L (i(w))$, then it is easy to see that
${\tilde E}_L$ can also be given by
\[
{\tilde E}_{L} = \lambda^\mu(C^\mu)(L_c) - L_c = \displaystyle
w^A\fpd{L_c}{w^A}(w) - L_c(w) \,\, \in \cinfty{\sW}.\]

The construction of $\theta_L$ in (\ref{thetaL}) can, however, not
directly be translated to $\mu^\rho$. Indeed, although ${\tilde
d}L_c$ is a well-defined 1-form on $\Sec(\mu^\rho)$, there exist
no vertical endomorphism on $\Sec(\mu^\rho)$. On the other hand,
we can also start from $\delta L_c$, which is a 1-form on
$\Sec(\mu^\lambda)$, but its image under the vertical endomorphism
$S^\mu$ of $\mu^\lambda$ gives us a 1-form
\[
\theta_{L_c} = S^\mu(\delta L_c) = \fpd{L_c}{w^A}(w){\mathcal X}^A
= i^a_A\fpd{L}{\sv^a}(i(w)){\mathcal X}^A
\]
 on $\Sec(\mu^\lambda)$ and
not on $\Sec(\mu^\rho)$, as required. Having found no direct
construction on $\mu^\rho$, it seems appropriate to take one more
step backwards (w.r.t.\ the diagrams of the previous section). We
will use a suitable restriction of $\theta_L$ to $\mu^\rho$.
\begin{dfn}
The Poincar\'e-Cartan 1-form ${\tilde\theta}_L$ is defined as the
1-form $(T^\rho i)^*(\theta_L)$ on $\Sec(\mu^\rho)$, i.e.\
\[
{\tilde\theta}_L(w)(\sv, W) = \theta_L(i(w))(T^\rho i(\sv,W)),
\] or, locally, $\displaystyle {\tilde\theta}_L =
\fpd{L}{\sv^a}(i(w)){\tilde{\mathcal X}}^a$. The energy function
${\tilde E}_{L}$ is the restriction of $E_L$ to $\sW$.
\end{dfn}
From the coordinate expressions it is easy to see that
$\theta_{L_c} =I^*({\tilde\theta}_L)$, i.e.\
\[\theta_{L_c}(w_1)(w_2,W) =
{\tilde\theta}_{L}(w_1)(I(w_2,W)). \] We now have all the
ingredients for a coordinate free description of
(\ref{extension}).
\begin{dfn}
A Lagrangian system on $\mu^\lambda:T^\lambda\sW\rightarrow\sW$ is
a pseudo-\sode\, $\Gamma\in \Sec(\mu^\lambda)$ that solves the
equation
\begin{equation}
\label{dynamics} i_{\Gamma} \delta {\tilde\theta}_L  = -\delta
{\tilde E}_{L}. \end{equation}
\end{dfn}

 It can easily be checked
that (\ref{dynamics}) gives indeed the correct equations
(\ref{extension}). First, one can calculate that
\[
\delta{\tilde\theta}_L = \left( \lambda^i_A i^b_B \frac{\partial
L^2}{\partial x^i \partial \sv^b } + \lambda^i_A i^b_B
\fpd{i^c_C}{x^i} \frac{\partial L^2}{\partial \sv^b \partial \sv^c
}  - \frac{1}{2} D^c_{AB} \fpd{L}{\sv^c}\right) {\mathcal X}^A
\wedge {\mathcal X}^B -i^b_B i^a_A \frac{\partial L^2}{\partial
\sv^a
\partial \sv^b} {\mathcal X}^A \wedge {\mathcal W}^B \]
and
\[
\delta {\tilde E}_{L} = \left( w^C i^a_C \frac{\partial
L^2}{\partial x^i
\partial \sv^a } + w^B w^C i^b_C \fpd{i^a_B}{x^i} \frac{\partial
L^2}{\partial \sv^a \partial \sv^b } - \fpd{L}{x^i} \right)
\lambda^i_A {\mathcal X}^A + w^Ai^a_A i^b_B \frac{\partial
L^2}{\partial \sv^a \partial \sv^b } {\mathcal W^B}.
\]
The coefficients $f^A$ of $\Gamma$ should therefore satisfy
\begin{equation}
 \frac{\partial^2
L}{\partial \sv^a\partial x^i} i^a_B \lambda^i_A  w^A +
\frac{\partial^2 L}{\partial \sv^a\partial \sv^b} \big(  w^B w^C
i^a_A \lambda^i_C \fpd{i^b_B}{x^i}+ i^a_B i^b_A f^A \big)-
\fpd{L}{x^i}\lambda^i_B - \fpd{L}{\sv^c} D^c_{BA}w^A =0 .
\label{Lagr1}
\end{equation}
The dynamics are given by the equations of the integral curves of
$X=\lambda^\mu(\Gamma)=w^A\lambda^i_A\fpd{}{x^i} + f^A
\fpd{}{w^A}$, i.e.\ they are solutions of \begin{equation}\left\{
\begin{array}{l}
\dot{x}^i \, = \, {\lambda}^i_{A}(x) { w}^{A},
\\[1mm]\dot{w}^A \, = \, { f}^{A}(x,w).
\end{array} \label{pseudosodeeq}
\right.
\end{equation}
Along such solutions, we find that
\[
\frac{d}{dt}\Big(\fpd{L_c}{w^A}\Big) = \Big(\frac{\partial^2
L}{\partial \sv^a
\partial x^i}i^a_A + \fpd{L}{\sv^a}\fpd{i^a_A}{x^i}\Big)\lambda^i_B w^B +
\frac{\partial^2 L}{\partial \sv^a\partial \sv^b}i^b_B i^a_A f^B.
\]
and
\[
\fpd{L_c}{x^i} = \fpd{L}{x^i} + \fpd{L}{v^a}\fpd{i^a_A}{x^i}w^A.
\]
 After plugging this information into (\ref{Lagr1}),
(\ref{pseudosodeeq}) becomes exactly (\ref{extension})
\begin{equation}\label{Lagr2}
\left\{
\begin{array}{rcl}
\dot{x}^i &\!\!\!=\!\!\!& {\lambda}^i_{A}(x) {w}^{A}, \\[2mm] \displaystyle
\frac{d}{dt}\Big(\fpd{L_c}{w^A}\Big) &\!\!\!=\!\!\!& \displaystyle
 \lambda^i_A\fpd{L_c}{x^i} + w^B\Big(D^c_{BA}  - \lambda^i_A
\fpd{i^c_B}{x^i} + \lambda^i_B \fpd{i^c_A}{x^i}\Big)
\fpd{L}{\sv^c}.
\end{array}
\right.
\end{equation}

\section{Examples}\label{examples}
The main example of this section shows that nonholonomic
mechanical systems with symmetry admit a formulation in the above
framework after reduction. The first two examples are
straightforward.

{\bf 1. Lagrangian systems on Lie algebroids.} If $\mu: \sW
\rightarrow M$ is exactly the Lie algebroid $\tau: \sV \rightarrow
M$, then $i^a_b = \delta^a_b$ and $L_c=L$. The equations
(\ref{Lagr2}) are then exactly (\ref{LAeq}). Examples of such
systems can be found e.g.\ in rigid body dynamics (see e.g.\
\cite{Mart} for a worked out example about the heavy top), control
theory \cite{Mart2} or Lagrangian systems with symmetry on
principal fibre bundles \cite{CMR,Wein}.

{\bf 2. Lagrangian equations on Lie subalgebroids.} With the help
of (\ref{substructure}), it is easy to see that in this case,
$\displaystyle \Big(D^c_{BA}  - \lambda^i_A \fpd{i^c_B}{x^i} +
\lambda^i_B \fpd{i^c_A}{x^i}\Big) \fpd{L}{\sv^c} = {D}_{BA}^{C}
\fpd{L_c}{{w}^{C}}$. The equations (\ref{Lagr2}) then become
exactly the equations (\ref{LAeq}) for the constrained Lagrangian
$L_c$ on the Lie algebroid $\mu:\sW\rightarrow M$. This is also
clear from the expression (\ref{dynamics}). Indeed, in this case
$\mu^\lambda: T^\lambda \sW \rightarrow \sW$ inherits the Lie
algebroid structure form $\mu$ and has also an exterior derivative
$d^\mu$ which acts on functions and forms in the following way
\[
d^\mu x^i = \lambda^i_A {\mathcal X}^A, \quad d^\mu w^A =
{\mathcal W}^A, \quad d^\mu {{\mathcal X}}^C = -\frac{1}{2}
D^C_{AB} {\mathcal X}^A {\wedge \mathcal X}^B \quad
\mbox{and}\quad d^\mu {{\mathcal W}}^A = 0.
\]
Obviously $\delta E_{L_c} = d^\mu E_{L_c}$. Moreover, due to
(\ref{substructure}),
\begin{eqnarray*}
\delta{\tilde\theta}_L &=& \delta\big(\fpd{L}{\sv^a}\big) \wedge
I^*({\tilde{\mathcal X}}^a) +I^*\big(\fpd{L}{\sv^a}\big)
\delta({\tilde{\mathcal X}}^a) = d^\mu \big(\fpd{L}{\sv^a}\big)
\wedge i_A^a{\mathcal X}^A - \frac{1}{2} \fpd{L}{\sv^a} D^a_{AB}
{\mathcal X}^A \wedge {\mathcal X}^B\\ &=& d^\mu
\big(\fpd{L}{\sv^a}\big) \wedge i_A^a{\mathcal X}^A - \frac{1}{2}
\fpd{L}{\sv^a} i_C^a D^C_{AB} {\mathcal X}^A \wedge {\mathcal X}^B
\\ &=& d^\mu \big(\fpd{L}{\sv^a}\big) \wedge i_A^a{\mathcal X}^A -
\frac{1}{2} \fpd{L_c}{w^C} D^C_{AB} {\mathcal X}^A \wedge
{\mathcal X}^B = d^\mu\theta_{L_c}.
\end{eqnarray*}
Therefore (\ref{dynamics}) is indeed $i_\Gamma d^\mu\theta_{L_c} =
- dE_{L_c}$.

{\bf 3. Non-holonomic systems.} In the introduction we already
have defined the Lagrange-d'Alembert equations describing a
non-holonomic mechanical system. We will use similar notations.
Suppose $Q$ is the configuration space of a mechanical system that
is subject to some kinematic (linear) constraints. If the
constraint distribution is $m$-dimensional, then, due to the
regularity condition of $D$, we can express $m$ velocities ${\dot
s}^\alpha$ defined by ${\dot x^i} = ({\dot r}^I, {\dot
s}^\alpha)$, up to a renumbering, in terms of the others
\[
{\dot s}^\alpha = - A^\alpha_I(x) {\dot r}^I.
\]

Let $\iota: D \rightarrow TQ$, then it is given by $(s^\alpha,
r^I,{\dot r}^I)\mapsto (s^\alpha, r^I, {\dot r}^I, {\dot s}^\alpha
= - A^\alpha_I {\dot r}^I)$. Before continuing, we wish to make
the following remark: in some geometrical models that treat
non-holonomic systems it is sometimes further assumed that there
exists a bundle structure of the configuration space $Q$ over some
manifold $N$ such that $D$ is the horizontal distribution $HM
\subset TQ$ of a connection on $Q \rightarrow N$ (this can always
be done locally). This additional assumption is not necessary for
our purposes. The purpose of this section is to show that the
equations for non-holonomic mechanical systems that can be found
in e.g.~\cite{Bloch}, fit into the framework presented above.

We can take $\sV$ to be $TQ$, equipped with the natural Lie
algebroid structure: the anchor map $\rho$ is the identity and the
Lie algebroid bracket is the usual Lie bracket of vector fields.
Then $\tau^\rho$ is nothing but $TTQ \rightarrow TQ$. The
subbundle $\sW$ is precisely the distribution $D$, i.e. $i=\iota$.
Likewise, $\mu^\rho$ becomes simply $T\sW \rightarrow \sW$. Since
$\lambda = i$, the last prolongation bundle $\mu^i$ (on which the
dynamics will be defined!) is $T^i\sW = \sW \times_i T\sW
\rightarrow \sW$. The Poincar\'e-Cartan 1-form ${\tilde\theta}_L$
is in this case, the pullback of the usual Poincar\'e-Cartan form
$\theta_L = \displaystyle \fpd{L}{{\dot r}^I} dr^I + \fpd{L}{{\dot
s}^\alpha} d s^\alpha$ for $L$ by the injection $i$, i.e.\ it is a
a differential form on $\sW$ which formally looks similar, but
where the coefficients should be evaluated along the constraints.
The biggest difference between our approach and many others is
that the fundamental form $\delta{\tilde\theta}_L$ and the
dynamics $\Gamma$ are {\em not} a differential form or a vector
field on $\sW$, but a form and a section of a prolongation bundle.
The equations (\ref{Lagr2}) (with $\lambda^I_J = \delta^I_J$ and
$\lambda^\alpha_J = - A^\alpha_J$) are of course the required
Lagrange-d'Alembert equations
\begin{equation}\label{nonhol}
\left\{
\begin{array}{rcl}
{\dot s}^\alpha &\!\!\!=\!\!\!& - A^\alpha_I {\dot r}^{I},
\\[2mm] \displaystyle \frac{d}{dt}\Big(\fpd{L_c}{{\dot r}^I}\Big) &\!\!\!=\!\!\!&
\displaystyle \fpd{L_c}{r^I} - A^\alpha_I \fpd{L_c}{s^\alpha} -
{\dot r}^J B^\alpha_{IJ} \fpd{L}{{\dot s}^\alpha}.
\end{array}
\right.
\end{equation}
where $\displaystyle B^\alpha_{IJ} = \fpd{A^\alpha_I}{r^J} -
\fpd{A^\alpha_J}{r^I} + A^\beta_I \fpd{A^\alpha_J}{s^\beta} -
A^\beta_J \fpd{A^\alpha_I}{s^\beta}$. The standard way to obtain
the equations (\ref{nonhol}) goes by taking (\ref{nonholeq}) as
starting point, and following the same procedure that has lead us
to (\ref{extension}) in the introduction (see also equations
(5.2.7) in \cite{Bloch}).

{\bf 4. The reduction of non-holonomic systems with symmetry.} In
this example, we will call the configuration space of the
mechanical system $Q$ and the Lagrangian $l$. We further assume
that $l$ is invariant under the (tangent lift of the) action of a
Lie group $G$ and that $\pG:Q\rightarrow Q/G$ has the structure of
a principle fibre bundle. We will first introduce the Lie
algebroid structure of interest. It is assumed that the reader is
familiar with the natural constructions associated with a
principal fibre bundle \cite{koba}.

Let $\tilde\la = (Q\times \la)/G \rightarrow Q/G$ be the
associated Lie algebra bundle (for a definition see e.g.\
\cite{Mac}, it is a Lie algebroid structure with vanishing anchor
map and structure functions the structure constants of the Lie
algebra $\la$). Suppose that $A: TQ \rightarrow \la$ is a
principal connection on $\pG$ with horizontal distribution $H$.
Then $TQ=H \oplus {V}\pG$, where
${V}\pG=\{(V\pG)_q=T_q(\pG^{-1}([q]))\, | \, q\in Q\}$ is the
vertical distribution. It is shown e.g.\ in \cite{CMR,Mac} that
the connection induces a vector bundle isomorphism $\alpha_A$
between $TQ/G$ and $T(Q/G)\oplus \tilde\la$ by means of
\[
\alpha_A([v_q]) = T\pG(v_q) \oplus [q \cdot A(v_q)], \qquad v_q
\in T_qQ.
\]

Here, $[q\cdot\xi]$ stands for the equivalence class of
$(q,\xi)\in Q \times \la$, i.e. $(q,\xi)\sim
(qg,Ad_{g^{-1}}\cdot\xi)$ for all $g\in G$. Remark that
$\alpha_A(H/G)= T(Q/G)$ and $\alpha_A(V\pG/G)=\tilde\la$. In the
notations from the previous sections, $M=Q/G$. The Lie algebroid
is then defined by $\tau:\sV=TQ/G\simeq T(Q/G)\oplus \tilde\la
\rightarrow Q/G$, where the bracket is taken to be the Lie bracket
restricted to right invariant vector fields on $Q$. The anchor map
$\rho:\sV \rightarrow TM$ of the Lie algebroid is nothing but the
projection onto $T(Q/G)$. The Lie algebroid bracket can be given
by
\begin{equation} \label{bracket1}
[X_1\oplus\ss_1, X_2\oplus\ss_2 ] = [X_1,X_2] \oplus
\big(\nabla_{X_1} \ss_2 - \nabla_{X_2} \ss_1 -\omega(X_1,X_2) +
[\ss_1,\ss_2]\big)
\end{equation}
with $X_i \in \vectorfields{M}$, $\ss_i \in \Sec(\tilde\la)$.
$\nabla$ denotes the covariant derivative on $\tilde\la\rightarrow
M$ associated with the connection $A$ and $\omega$ is its
curvature.

Recall that a local trivialisation of the principal fibre bundle
$\pG$ induces a bundle adapted coordinate chart on every
associated bundle of $\pi_G$. In particular, if $e_\alpha$ is a
basis for the Lie algebra $\la$, then a basis for the set of local
sections of $\tilde\la$ are defined by ${\ov e}_\alpha(x) =
[\psi^{-1}(x,e)\cdot e_\alpha]$ where $x\in U\subset M=Q/G$ and
$\psi: \pi_G^{-1}(U) \to U\times G$ is such a local trivialisation
of $\pi_G$. Suppose $(x^i, {\ov v}^\alpha)$ are coordinates on
$\tilde\la$ w.r.t.\ this basis. Then, $\sV$ has coordinates
$(x^i,\sv^a = ({\dot x}^i,{\ov v}^\alpha))$. Using this coordinate
system, the coefficients of the Lie bracket $[\ss_1,\ss_2]$ on
$\tilde\la$ are exactly the structure constants
$C^\gamma_{\alpha\beta}$ of $\la$. Furthermore, the connection
coefficients of the covariant derivative $\nabla$ on
$\tilde\la\rightarrow Q/G$ associated with $A$ take the form
$\Gamma_{i\alpha}^\beta = C^\beta_{\alpha\gamma}A^\gamma_i$.
Finally, the coefficients of the $\tilde\la$-valued 2-form
$\omega$ on $M$ are
$$\displaystyle{\omega^\alpha_{ij} = \fpd{A^\alpha_j}{x^i} -
\fpd{A^\alpha_i}{x^j} + C^\alpha_{\beta\gamma} A^\beta_j
A^\gamma_i}.$$

The basis $\{{\ov e}_\alpha\}$ for $\Sec(\tilde\la)$ induces a
basis $\{e_i = \fpd{}{x^i}\oplus 0, e_\alpha = 0 \oplus {\ov
e}_\alpha \}$ for $\Sec(\tau)$. W.r.t.\ this basis, the Lie
algebroid bracket (\ref{bracket1}) is given by
\[
[e_i,e_j] = - \omega^\gamma_{ij} e_\gamma, \qquad [e_i,e_\alpha] =
\Gamma^\gamma_{i\alpha}e_\gamma, \qquad [e_\alpha,e_\beta]  =
C^\gamma_{\alpha\beta} e_\gamma.
\]
This Lie algebroid structure on $\tau:TQ/G\simeq TM\oplus
\tilde\la\rightarrow M$ is the so-called {\em Atiyah algebroid}
(see also \cite{DMM,Mac}).

Suppose now that the system is subject to some linear constraints
$D\subset TQ$. In contrast with the previous example, there is a
natural fibration of $Q$ available, $Q \rightarrow M=Q/G$, so it
makes sense to compare $D$ with the vertical subspace ${V}\pG$ of
this fibration, rather than to assume that $D$ is the horizontal
distribution of a connection on some fibre bundle. We will follow
here the approach of \cite{CMR2}, although there are many other
\cite{Bloch2,Frans1,Frans2,Cortes}. In \cite{CMR2} two additional
assumptions were made. First, we suppose that $T_qQ=D_q +
(V\pG)_q$. This assumption means, among others, that $S=\{S_q= D_q
\cap (V\pG)_q\, | \, q\in Q\}$ is a subbundle of $TQ$, $D$ and
$V\pG$. Further, we will also assume that $D$ (and therefore also
$S$) is $G$-invariant. In \cite{CMR} it is proved that there
exists always a $G$-invariant metric on $Q$, which we now assume
to be fixed. Let $H_q$ be the orthogonal complement of $S_q$ in
$D_q$ (with respect to this metric), then $D=H \oplus S$, while
$TQ = H\oplus V\pG$. Further, if $U$ is the orthogonal complement
of $S$ in $V\pG$, then $TQ=H\oplus S\oplus U$. Due to the above
assumptions, all three distributions are $G$-invariant and thus
$TQ/G=H/G\oplus S/G\oplus U/G$. Let $A:TQ\rightarrow \la$ be the
principal connection whose horizontal subspace at $q$ is exactly
$H_q$. We now use this connection to consider the above
decomposition of $TQ/G$ in the isomorphic bundle $TM\oplus \tilde
\la$, i.e.\ we have \begin{eqnarray*} && TM\oplus 0
=\alpha_A(H/G),\\ && 0 \oplus
\tilde\la=\alpha_A(V\pG/G)=\alpha_A(S/G)\oplus \alpha_A(U/G) =
\tilde \sa \oplus \tilde \ua,\\ && \alpha_A(D/G)=TM\oplus
\sa.\end{eqnarray*}

Analogously to the reduction of non-constrained systems with
symmetry, the Lagrange-d'Alembert principle which describes the
equations of motion on $Q$, can be reduced to a `variational
principle' on a reduced space (cf. \cite{CMR2}). This new
principle generates the reduced equations, the so-called
Lagrange-d'Alembert-Poincar\'e equations. We now show that these
equations are Lagrange equations on a subbundle of a Lie
algebroid. Using the above notations we now define $\mu:\sW=
D/G\simeq TM \oplus \tilde \sa \rightarrow M$ and this bundle is
the required subbundle of the Atiyah algebroid $\tau:\sV = TQ/G
\simeq TM\oplus \tilde \la \rightarrow M$. Choose a basis $\{{\ov
e}_I\}$ of $\Sec(\tilde\sa)$ and denote ${\ov e}_I = e^\alpha_I
{\ov e}_\alpha$. Let $(x^i,{\ov w}^I)$ be coordinates on $\sa$,
then $i:\sW \rightarrow\sV$ has components $i_k^j =\delta_k^j$,
$i^j_I=0$, $i^\alpha_i =0$ and $i_J^\alpha = e^\alpha_J$.
Moreover, the components for $\lambda:\sW\rightarrow TM$ are
$\lambda^i_j=\delta^i_j$, $\lambda^i_J = 0$.

Let $L\in\cinfty{\sV}$ be the reduced Lagrangian, i.e.\
$L([v_q])=l(v_q)$. If $L_c$ is the restriction of $L$ to $D/G$,
then the Lagrange-d'Alembert-Poincar\'e equations are
\begin{equation} \label{LPeq} \left\{
\begin{array}{rll} \displaystyle{ \frac{d}{dt}\fpd{L_c}{{\ov w}^I}} &=&
\displaystyle{ -\fpd{L}{{\ov v}^\beta}\Big(D^\beta_{IJ}{\ov w}^J -
\Gamma^\beta_{j I}{\dot x}^j + \fpd{e_I^\beta}{x^j} {\dot x}^j
\Big),} \\[3mm] \displaystyle{\frac{d}{dt} \fpd{L_c}{{\dot
x}^i} - \fpd{L_c}{x^i} } &=& \displaystyle{-\fpd{L}{{\ov
v}^\beta}\Big(\Gamma^\beta_{iJ}{\ov w}^J - \omega^b_{ij}{\dot x}^j
- \fpd{e^\beta_J}{x^i}{\ov w}^J \Big).}
\end{array}\right.
\end{equation}
with $D^\beta_{IJ} = C^\beta_{\gamma\delta} e^\gamma_I e^\gamma_J$
and $\Gamma^\beta_{jI}= \Gamma^\beta_{j\alpha} e^\alpha_I$ (see
also equations (5.8.45-47) in \cite{Bloch} and (4.11-4.14) in
\cite{CMR2}). In \cite{Bloch, Frans2, CMR2,Cortes} one can find
examples considered in full detail of nonholonomic systems with
symmetry, for instance the snakeboard and the vertically rolling
disk.

{\bf 5. Normal extremals in Lagrangian systems on Lie algebroids.}
The last example can be found in the theory of geometric optimal
control theory. We first briefly recall some general concepts from
control theory. In control theory one studies dynamical systems
that can be steered by external devices (typically representing a
human input to the system). In dynamical systems theory, such
systems are typically represented by a differential equation of
the following type
$$\dot x(t) =f(x(t),u(t)),$$ where $x\in \R^n$ represents the
configuration of the system and where $u \in \R^k$ represents the
control functions. It should be clear that the dynamics of the
system is completely determined (up to the initial condition) by
the control $u(t)$. A geometric framework for studying control
theory is that of an anchored bundle, i.e. we assume that the
configuration space is a manifold $M$ and that the control space
$\R^k$ is the $k$-dimensional fibre of a bundle $\sV$ over $M$,
with projection $\tau$. The analogue of the map $f$ is a bundle
map $\rho$ from $\sV$ to $TM$, fibred over the identity. An
admissible curve is a curve $\sv(t)$ in $\sV$ such that $d/dt
(\tau(\sv(t))) = \rho(\sv(t))$. It is not difficult to see that in
a local coordinate system this condition has the precise structure
a control differential equation. Assume that a Lagrangian function
$L \in \cinfty{\sV}$ is given. In optimal control theory one
studies the following variational problem: `among all admissible
curves $\sv$ defined over the interval $[t_0,t_1]$ and such that
$\tau(\sv(t_0))=x_0$ and $\tau(\sv(t_1))=x_1$ for fixed endpoints
$x_0,x_1\in M$, which one minimises a cost functional $\int
L(\sv(t))dt$~?'. The maximum principle \cite{pont} gives necessary
conditions for admissible curves to be minimising. Locally they
are given as follows: an admissible curve
$\sv(t)=(x^i(t),\sv^a(t))$ is minimising if there is (1) a curve
$p_i(t)$ in $T^*M$ along $x^i(t)$ and a constant real number
$p_0=0,-1$ such that $(p_0,p_1(t),\ldots,p_n(t))\neq 0 \ \ \forall
t$ and (2) the following Hamiltonian system is satisfied at all
time $t$, with $H(x^i,\sv^a,p_i)=p_i\rho^i(x^j,\sv^a) + p_0
L(x^i,\sv^a)$ a function on $T^*M\times\sV$:
\begin{eqnarray*}
\dot x^i(t) &=& \fpd{H}{p_i} = \rho^i(x^j(t),\sv^a(t)),\\
\dot p_i(t)&=& -\fpd{H}{x^i} = -p_j\fpd{\rho^j}{x^i} -p_0 \fpd{L}{x^i}\\
0 &=& \fpd{H}{\sv^a}=p_i\fpd{\rho^i}{\sv^a} +p_0 \fpd{L}{\sv^a}.
\end{eqnarray*}
The latter condition says that the function $\sv^a\mapsto
H(x^i(t),p_i(t),\sv^a)$ on the fibres of $\sV$ attains a local
maximum at the point $\sv^a=\sv^a(t)$. A coordinate free version
for this theorem was proven by H.J.\ Sussmann in \cite{sussmann4}.
However, for the sake of simplicity, we continue to work in a
local coordinate system. It is very interesting to note that there
are two kinds of solutions: those admissible curves for which
there is a $p_i(t)$ satisfying the conditions with $p_0 < 0$ and
those for which $p_0 =0$. The latter are called abnormal extremals
since in this case the conditions from the maximum principle do
not depend on the cost function (see also
\cite{control,Mont,Montboek}). If $p_0< 0$, then the solutions are
called normal extremals. We now assume that the bundle $\sV$ has
the structure of a Lie algebroid, and that the map $\rho$ is the
anchor map. We will show that the above conditions from the
maximum principle for normal admissible curves can be rewritten as
the solutions to Lagrangian systems on Lie algebroids with
constraints. For that purpose we consider the equation expressing
the maximality condition:
\[
0=\fpd{H}{\sv^a}(x^i(t),\sv^a(t),p_i(t)) = p_i(t)\rho^i_a(x^j(t))
- \fpd{L}{\sv^a}(x^j(t),\sv^b(t)).
\]

Since we assumed that the Lagrangian is regular, we can consider
the inverse of $\F L= \partial L/\partial \sv^a: \sV \to \sV^*$,
which gives us the following condition on the control curve:
\begin{equation}\label{vgllegendre}
\sv^a(t) = (\F L^{-1})^a(x^j(t),p_i(t)\rho^i_a(x^j(t))).
\end{equation}

If we substitute this in the Hamiltonian equations from the
maximum principle, we obtain after some straightforward
calculations (and taking into account the structure equations of
the Lie algebroid $\sV$) that $(x^i(t),\sv^a(t))$ is a solution to
the Lagrangian equations on the Lie algebroid:
\begin{equation} \left\{
\begin{array}{rcl}
\dot{x}^i &\!\!\!=\!\!\!& {\rho}^i_{a}(x) {\sv}^{a}, \\[2mm] \displaystyle
\frac{d}{dt}\left(\fpd{L}{{\sv}^{a}}\right) &\!\!\!=\!\!\!&
\displaystyle {\rho}^i_{a} \fpd{L}{x^i} - {C}_{ab}^{c} {\sv}^{b}
\fpd{L}{{\sv}^{c}}.
\end{array}
\right.\label{vgllaatste}
\end{equation}
This system is in general not equivalent with the Hamiltonian
system from the maximum principle. The condition that $\sv^a$ is
in the image of $\F L^{-1}$ has to be taken into account. In the
specific case that $\F L^{-1}$ is a linear map, this is precisely
saying that $\sv(t)$ is contained in a linear subbundle $\sW$,
i.e. $\sv(t) \in \sW =\mbox{im}(\F L^{-1}\circ i(T^*M))$.
Solutions of the systems (\ref{vgllegendre}) and
(\ref{vgllaatste}) are in fact special solutions of the
system~\eqref{NHeq}. So, normal extremals are solutions of a
Lagrangian system on a subbundle of a Lie algebroid. However, the
equations of motion for the normal extremals are {\em not
equivalent} to the non-holonomic equations of motion. The normal
extremals have to satisfy stronger conditions \eqref{vgllaatste}.

A sufficient condition for the image of $\F L^{-1}$ to be a linear
subbundle $\sW$ of $\sV$ is the condition that $L$ is a Lagrangian
of mechanical type, i.e.\ when $L=T -\tau^*V$, where $T$ is the
kinetic energy associated with a metric on $\sV$ and where $V$ is
a potential function defined on $M$.  A typical example of a
mechanical system on a Lie algebroid is the spinning top (cf.
\cite{Mart}).

\section{Conclusions and Outlook}

In this contribution we have dealt with systems on Lie algebroids
that are subject to some constraints. We obtained an intrinsic
description of the dynamics of the systems in terms of a section
of a prolongation bundle. Our approach unifies models for both
Lagrange-d'Alembert equations and Lagrange-d'Alembert-Poincar\'e
equations.

Such a simultaneous description can be very handy in many
applications. For example, in \cite{CMR} it is shown that the
so-called `Lagrange-Poincar\'e bundles' form the ideal platform on
which the reduction process of (unconstrained) Lagrange-Poincar\'e
equations can be repeated. Knowing that such bundles are in fact
Lie algebroids, we hope that Lie algebroid theory (and in
particular the description of such systems as sections of a
prolongation Lie algebroid) will play an important role in future
developments of the process called `Lagrangian reduction by
stages'. Due to the observations in this paper, it has become
clear that, if we want to find a geometric formalism for
successive reduction of non-holonomic systems, we will need to
explore the geometry of (the prolongation of) Lagrange-Poincar\'e
subbundles.

Also in the case of non-holonomic systems (example~3), the above
developped theory leads to interesting new insights. Indeed, the
prolongation bundle $T^i\sW \rightarrow \sW$, where e.g. the
Poincar\'e-Cartan 2-form $\delta{\tilde\theta}_L$ and the dynamics
$\Gamma$ live, plays a crucial role in our approach, but its
importance has, so far, not been recognized in the literature. A
framework that seems to be closely related to ours is
\cite{Willy1,David} (although their set-up is more general since
also time-dependent systems were included). In those papers, two
fundamental two-forms on $T\sW$ have been considered. It would be
of interest to explore the relations between those two forms on
the one hand and $\delta{\tilde\theta}_L$ on the other hand.
Further, we intend to find out how our operator $\delta$ fits in
the approach of \cite{Willy1,David}.

An other path for future developments is that of an appropriate
framework for studying Hamiltonian equations on a subbundle of a
Lie algebroid. Hamiltonian systems on Lie algebroids were already
considered in e.g. \cite{DMM, Mart2}. If $\tau^*:\sV\rightarrow M$
and $\mu^*:\sW^*\rightarrow M$ are the duals of $\tau$ and $\mu$,
then the main object of a Hamiltonian description for constraint
systems will be a section of the prolongation bundle
$(\mu^*)^\lambda:T^\lambda\sW^*\rightarrow\sW^*$. Similar as
before, also the, not unrelated, bundles $(\tau^*)^\rho$,
$(\tau^*)^\lambda$ and $(\mu^*)^\rho$ will come into the picture.
In the special example of systems with symmetry, we should be able
to relate the first equation in (\ref{LPeq}), in a Hamiltonian
formulation, with the momentum equation (for a recent survey see
\cite{Bloch}).

{\bf Acknowledgements.} We are indebted to Frans Cantrijn, Eduardo
Mart\' inez and Willy Sarlet for useful discussions.

\end{document}